\documentclass[12pt]{article}
\usepackage{amssymb,amsmath}
\def\qed{\hfill$\Box$\par}
\def\a{\alpha}

\def\D{\Delta}
\def\UU{{\cal U}}
\def\FF{{\cal F}}
\def\LL{{\cal L}}

\def\ad{{\rm ad\,}}
\def\Aut{{\rm Aut}}

\def\exp{{\rm exp}}

\def\e{\epsilon}

\def\ssc{\scriptscriptstyle}

\def\cl{\centerline}

\def\vs{\vspace*}

\def\UU{{\cal U}}
\def\WW{{\cal W}}
\def\LL{{\cal L}}
\def\FF{{\cal F}}

\def\ni{\noindent}

\def\N{\mathbb{N}{\ssc\,}}
\def\Z{\mathbb{Z}{\ssc\,}}

\def\F{\mathbb{F}{\ssc\,}}
\def\QED{\hfill$\Box$} \numberwithin{equation}{section}
\newtheorem{theo}{Theorem}[section]
\newtheorem{defi}[theo]{Definition}

\newtheorem{lemm}[theo]{Lemma}

\def\adddot{$\!\!\!${\bf.}\ \ }

\begin{document}

\cl{{\large \bf On Drinfel'd twist deformation of the super-Virasoro
algebra}\footnote{ Supported by Science $\&$ Technology Program of
Shanghai Maritime University.\\[2pt] \indent Corresponding E-mail:
hyyang@shmtu.edu.cn}} \vs{6pt}

\cl{ Hengyun Yang}

\cl{\small Department of Mathematics, Shanghai Maritime University,
Shanghai 201306, \vs{-4pt}China }

\vs{6pt}

\noindent{\small{\bf Abstract.} In this paper, we describe
nonstandard quantum deformation of the super-Virasoro algebra. Using
the Drinfel'd twist quantization technique, we obtain the deformed
coproduct and antipode. Hence we get a family of noncommutative and
noncocommutative Hopf superalgebras.

  }


\noindent{\bf1. \
Introduction}\setcounter{section}{1}\setcounter{equation}{0}

Quantum groups, mathematically carry the structures of
noncommutative and noncocommutative Hopf algebras, were first
introduced by Drinfel'd and Jimbo. One of the most important
examples of quantum groups is deformation of the universal
enveloping algebra of a Lie algebra. This deformation induces a Lie
bialgebra structure on the underlying Lie algebra. In \cite{D},
Drinfel'd posed the quantization problem of the Lie bialgebra.
Lately, Etingof and Kazhdan \cite{EK} settled this question. But a
general formula of a quantization was not obtained. Many authors
have made great efforts to quantize explicitly some Lie bialgebras.

Inspired by the discovery of quantum groups, quantum supergroups,
i.e., Hopf superalgebras, have also been defined (c.f. \cite{BGZ,
CK}), which provide a powerful tool for constructing trigonometric
solutions of the $\Z_2$-graded Yang-Baxter equation. By extending
Etingof and Kazhdan's work \cite{EK}, Geer \cite{G1} recently proved
that there exist a general quantizations of Lie superbialgebras.
Similar to Lie algebras cases, the deformations of Lie superalgebras
are not unique. In \cite{Z}, Zhang proved that there is a new Hopf
superalgebra structure by the Drinfel'd twist. Using this method,
some good Hopf superalgebras have been found in recent years, e.g.,
Aizawa \cite{Ai} and Celeghini et al \cite{CK1} studied the
drinfel'd twist deformations of $sl(1|2)$ and $osp(1|2)$,
respectively. By the Drinfel'd twist, \cite{G} and \cite{YS} give
two different Hopf algebra structures on the Witt algebra. The aim
of this paper is to construct the quantization of the super-Virasoro
algebra, which is generated by the same Drinfel'd twist of the Witt
algebra, studied in \cite{Y} ( see also \cite{SZ, S, Y}). As a
by-product, we obtain two combinational identities (see
(\ref{combin1}) and (\ref{combin2})).

Throughout this paper, $\F$ denotes a field of characteristic zero,
all vector space and tensor products are over $\F$.  Let $\Z$,
$\Z_+$, $\N$ denote the sets of all integers, nonnegative integers,
positive integers, respectively. We use the convention that if an
undefined term appears in an expression, we always treat it as zero;
for instance, $L_{\frac12}=0$ if $\frac12\notin\Z$.

\vskip12pt

\noindent{\bf2. Main results
}\setcounter{section}{2}\setcounter{theo}{0}\setcounter{equation}{0}

Now let us start by recalling some definitions and preliminary
results. A supervector space $H$ is a $\Z_2$-graded vector space,
i.e., $H=H_{\bar 0}\oplus H_{\bar 1}$. If an element $x$ is in
either $H_0$ or $H_1$, we say that it is $\Z_2$-homogeneous. We
assume that all elements below are $\Z_2$-homogeneous, where
$\Z_2=\{\bar 0, \bar 1\}$. For $x\in H$, we always denote
$[x]\in\Z_2$ to be its parity, i.e., $x\in H_{[x]}$. We say that $x$
is even (odd) if $x\in H_{\bar 0}$ (resp. $x\in H_{\bar 1}$). A
superalgebra $(H,\mu,\tau)$ over a commutative ring $R$ is a
supervector space equipped with an associative product $\mu
:H\otimes H\rightarrow H$ respecting the grading and a unit element
$1\in H_{\bar 0}$. A Hopf superalgebra $(H, \mu, \tau, \D, \e, S)$
is a superalgebra equipped with a coproduct $\D:H\rightarrow
H\otimes H$, a counit $\e: H\rightarrow \F$, and an antipode
$S:H\rightarrow H$, satisfying compatibility conditions. Note that
$S$ satisfies $S(xy)=(-1)^{[x][y]}S(y)S(x)$ for $x,y\in H$.
\begin{defi} A Drinfel'd twist $\cal{F}$ is an invertible element of $ H \otimes H$
and satisfies
\begin{eqnarray}\label{2.10}
&\!\!\!\!\!\!\!\! &(\FF \otimes 1) (\D \otimes Id ) (\FF ) = (1
\otimes \FF) (1 \otimes \D ) (\FF ),\\
 \label{2.11}&\!\!\!\!\!\!\!\! & (\e \otimes Id )
(\FF ) = 1 \otimes 1 = (Id \otimes \e ) (\FF) .\end{eqnarray}
\end{defi}

Write
$$\FF=\sum\limits f_{(1)}\otimes f_{(2)},\quad \FF^{-1}=\sum\limits f'_{(1)}\otimes f'_{(2)},$$
and set
$$u= \mu \cdot (S \otimes Id ) (\FF^{-1})=\sum\limits
S(f'_{(1)})f'_{(2)}.$$ The property of $S$ shows that $u$ is
invertible with inverse
$$u^{-1}=\mu \cdot (Id \otimes S )(\FF)=\sum\limits
f_{(1)}S(f_{(2)}).$$

The following result gives a method to construct new Hopf
superalgebra from old ones (cf.~\cite{Z}), the non-super case can be
found in \cite{D}.

\begin{lemm}\label{lemma2.3} The superalgebra $(H, \mu, \tau, \tilde{\D}, \tilde{\e}, \tilde {S} )$ is a new Hopf
superalgebra with
$$\tilde{\D} = \FF \D \FF^{-1} , \ \ \ \
\,\tilde{\e}=\e,\ \ \ \ \   \tilde{ S} = u^{-1} S u.$$
\end{lemm}

Now let us recall that the classical super-Virasoro algebra $\LL$
without central extension over $\F$ is defined as an
infinite-dimensional Lie superalgebra generated by the generators
$\{L_i,G_k\,|\,i\in\Z, k\in \frac12\Z\}$ satisfying the defining
relations
\begin{equation*}
\label{e-SVir} [L_i,L_j]=(j-i)L_{i+j},\quad
[L_i,G_k]=(k-\frac{i}{2})G_{i+k},\quad [G_k,G_l]=2L_{k+l},
 \end{equation*} for all
$i,j\in\Z,\,k,l\in \frac12 Z.$ Obviously, $\LL$ contains the Witt
algabra $\WW$ as subalgebra. In the following, we fix
$m\in\Z/\{0\}$, $\a\in\F$. Denote
$$X=\frac{1}{m} (L_0+\a m L_{-m}),\ \ \ Y={\rm
exp}(\a{\rm ad} L_{-m})(L_{m}).$$ In \cite{SZ}, Su and Zhao proved
that $X$ and $Y$ span a two-dimensional subalgebra of the Virasoro
algebra, i.e., $[X,Y]=Y$.

Denote by $\UU (\LL)$ the universal enveloping algebra of $\LL$. For
any $x \in \UU (\LL), a\in \F,r,k\in\Z_+, i\in\Z,$ we set
\begin{eqnarray*}
\label{2.1}&\!\!\!\!\!\!\!\! &x_a ^{<r>} = (x + a)(x+ a+1)\cdots
(x+a + r-1),\quad x_a ^{[r]} = (x+ a)(x+a-1)\cdots (x+a-r+1),
\\[6pt]
\label{2.3}&\!\!\!\!\!\!\!\! &\binom{a}{r}=\frac{a(a-1)\cdots(a-r+1)
}{r!},\qquad\qquad\qquad\quad \left[\begin{array}{c} a \\ r
\end{array} \right]_k=\frac{a(a-k)(a-2k)\cdots \big(a-(r-1)k\big)}{r!}.
\end{eqnarray*}
Especially, we have $\left[\begin{array}{c} a \\ r
\end{array} \right]_1=\left(\begin{array}{c} a \\ r
\end{array} \right) $ and $\left[\begin{array}{c} a \\ r
\end{array} \right]_{-1}=\left(\begin{array}{c} 1-a-r \\ r
\end{array} \right)$. Denote $x^{<r>} = x_0 ^{<r>} ,
 x^{[r]} = x_0 ^{[r]}$.
Note that $\UU (\LL)$ has a natural $\Z_2$-graded Hopf superalgebra
structure. On the generators $x\in \LL$ and the unit element $1$, we
define
\begin{equation*}\begin{array}{ll}
\D_0 (x ) = x\otimes 1 + 1 \otimes x , & \D_0(1)=1\otimes 1,\\
\e_0(x) =0, & \e_0(1)=1,\\
S_0(x) = -x & S_0(1)=1.
\end{array}\end{equation*}
Obviously, the Hopf superalgebra
$(\UU(\LL),\mu,\tau,\Delta_0,S_0,\epsilon_0)$ is cocommutative.
Denote by $\UU (\LL)[[t]]$ an associative $\mathbb{F}$-algebra of
formal power series with coefficients in $\UU (\LL)$. Then $\UU
(\LL)[[t]]/t\UU (\LL)[[t]]\\ \cong \UU (\LL)$. Naturally, $\UU
(\LL)[[t]]$ equips with an induced Hopf superalgebra structure
arising from that on $\UU (\LL)$, denoted still by $(\UU
(\LL)[[t]],\mu,\tau,\Delta_0,S_0,\epsilon_0)$, called the quantized
universal enveloping superalgebra. As the non-super case \cite{D},
Andruskiewitsch \cite{A} proved that the Lie superalgebra $\LL$ has
a natural Lie superbialgebra structure. Thus $\UU (\LL)[[t]]$ is
also called the quantization of the Lie superbialgebra $\LL$.
Define the coproduct $\D$ and the antipode $S$ on $\UU (\LL)[[t]]$
as follows:
\begin{eqnarray}
\D(L_i)\!\!\!& =\!\!\!&
\sum\limits_{r=0}^{\infty}\a^r\left[\begin{array}{c} (r-2)m-i \\
r
\end{array} \right]_m\big(L_{i-rm}\otimes (1-Yt)^{\frac{i}{m}-r}Y^rt^r\big)\label{main1}\\
\!\!\!& \!\!\!&
+\sum\limits_{r=0}^{\infty}(-1)^r\sum\limits_{s=0}^{2r}\a^s
a_s(r,i)\big(X^{<r>}\otimes
(1-Yt)^{-r}L_{i+(r-s)m}t^r\big),\nonumber \\
\D(G_k)\!\!\!& =\!\!\!&\sum\limits_{r=0}^{\infty}\a^r\left[\begin{array}{c} (r-\frac32)m-k \\
r
\end{array} \right]_m\big(G_{k-rm}\otimes (1-Yt)^{\frac{k}{m}-r}Y^rt^r\big)\label{main2}\\
\!\!\!& \!\!\!&
+\sum\limits_{r=0}^{\infty}(-1)^r\sum\limits_{s=0}^{2r}\a^s
b_s(r,k)\big(X^{<r>}\otimes
(1-Yt)^{-r}G_{k+(r-s)m}t^r\big),\nonumber\\ \e(L_i)\!\!\!& =\!\!\!&
0,\ \ \ \ \ \ \ \ \ \e(G_k)=0,\label{mian3}
\end{eqnarray}
\begin{eqnarray} S (L_{i} )\!\!\!& =\!\!\!& -(1-Yt)^{-\frac{i}{m}}\sum\limits_{r=0} ^{\infty}
\sum\limits_{p=0} ^{\infty}\sum\limits_{q=0}^{2p}(-1)^{r}\a^{r+q}\left[\begin{array}{c} (r-2)m-i \\
r
\end{array} \right]_m\times\label{main4}\\
\!\!\!& \!\!\!&\times
a_q(p,i-rm)X_{-\frac{i}{m}}^{[p]}L_{i+(p-r-q)m}Y ^{r} t^{r+p},\nonumber\\
 S (G_k )\!\!\!& =\!\!\!&-(1-Yt)^{-\frac{k}{m}}\sum\limits_{r=0} ^{\infty}
\sum\limits_{p=0} ^{\infty}\sum\limits_{q=0}^{2p}(-1)^{r}\a^{r+q}\left[\begin{array}{c} (r-\frac32)m-k \\
r
\end{array} \right]_m\times\label{main5}\\
\!\!\!& \!\!\!&\times
b_q(p,k-rm)X_{-\frac{k}{m}}^{[p]}G_{k+(p-r-q)m}Y ^{r}
t^{r+p},\nonumber
\end{eqnarray}
where $i\in\Z$, $k\in\frac12\Z$, and for $r\in\Z_+$, $0\leq s\leq
2r$,
\begin{equation}
\label{2.3'} a_s(r,i)=
 \sum\limits_{p=0}^{s}(-1)^p
\left[\begin{array}{c} i+m \\ p
\end{array}
\right]_m \left[\begin{array}{c} i+(r-p-2)m \\ r
\end{array}
\right]_{m}\left[\begin{array}{c} i+(r-p+1)m \\ s-p
\end{array}
\right]_{m},
\end{equation}
\begin{equation}\label{2.3''} b_s(r,k)=\sum\limits_{p=0}^{s}(-1)^p\left[\begin{array}{c} k+\frac{m}{2} \\ p
\end{array}
\right]_m \left[\begin{array}{c} k+(r-p-\frac32)m \\ r
\end{array}
\right]_{m} \left[\begin{array}{c} k+(r-p+\frac12)m \\ s-p
\end{array}
\right]_{m}.\end{equation}

The main result of this paper is the following theorem which gives
the quantization of the super-Virasoro algebra.
\begin{theo}\adddot\label{main}
The superalgebra $\UU (\LL )[[t]]$ under the the coproduct $\D$, the
counit $\e$, and the antipode $S$ defined by
(\ref{main1}-\ref{main5}) is a noncommutative and noncocommutative
Hopf superalgebra.
\end{theo}

As a by-product of our proof, we obtain the following two
combinatorial identities, which do not seem to be obvious to us:
\begin{eqnarray}
\label{combin1} \!\!\!&\!\!\!&\sum\limits_{p=0}^{s}(-1)^p
\left[\begin{array}{c} i+m \\ p
\end{array}
\right]_m \left[\begin{array}{c} i+(r-p-2)m \\ r
\end{array}
\right]_{m}\left[\begin{array}{c} i+(r-p+1)m \\ s-p
\end{array}
\right]_{m}=0, \\
\label{combin2}\!\!\!&
\!\!\!&\sum\limits_{p=0}^{s}(-1)^p\left[\begin{array}{c}
k+\frac{m}{2} \\ p
\end{array}
\right]_m \left[\begin{array}{c} k+(r-p-\frac32)m \\ r
\end{array}
\right]_{m} \left[\begin{array}{c} k+(r-p+\frac12)m \\ s-p
\end{array}
\right]_{m}=0,\end{eqnarray} where $i\in\Z,k\in\frac12\Z,
r\in\Z_+,\,s>2r.$

\vskip7pt

\noindent{\bf3. Proof of Theorem
\ref{main}}\setcounter{section}{3}\setcounter{theo}{0}\setcounter{equation}{0}

We shall divide the proof of Theorem \ref{main} into several lemmas.
Set
$$X'=\frac{1}{m}L_0,\ \ \  Y'=L_{m}.$$ Let
$$\overline{\WW}=\{\sum_{i\in\Z} a_iL_i |\, a_i\in\F, \mbox{ and }
a_i=0 \mbox{ for } i\ll 0 \} $$ be the {\it completed Witt Lie
algebra}. Then $\exp {(\ad L_{-m})}\in \Aut (\overline{\WW})$.
Evidently, we have
$$\exp{(\a\ad L_{-m})}(X')=X, \ \ \ \ \ \exp{(\a\ad L_{-m})}(Y')=Y.$$

\begin{lemm}\label{lemm3.1} For any  $r\in\Z_+, i\in\Z,$ we have
\begin{eqnarray}\label{3.4'}(\ad Y)^r(L_i) = \sum\limits_{q=0}^{2r} \a^q r!a_q(r,i)
L_{i+(r-q)m},\\
\label{3.4''}(\ad Y)^r(G_k) = \sum\limits_{q=0}^{2r} \a^q r!b_q(r,i)
G_{k+(r-q)m}.
\end{eqnarray}
 \end{lemm}
\ni{\it Proof.~} Note that $\exp {(\a\ad L_{-m})}\in \Aut
(\overline{\WW})$ with the inverse $\exp {(-\a\ad L_{-m})}$. We have
\begin{eqnarray}
(\ad Y)^r(L_i) \!\!\!&= \!\!\!&\exp {(\a\ad L_{-m})} \exp {(-\a\ad L_{-m})}((\ad Y)^r(L_i))\nonumber\\
\!\!\!&= \!\!\!&\exp {(\a\ad L_{-m})}(\ad Y')^r\exp {(-\a\ad L_{-m})}(L_i)\nonumber\\
\!\!\!&= \!\!\!&\exp {(\a\ad L_{-m})}(\ad Y')^r\sum\limits_{p=0}
^{\infty} (-\a)^{p}\left[\begin{array}{c} i+m \\ p
\end{array}
\right]_m L_{i-pm}\nonumber\\
\!\!\!&= \!\!\!&\exp {(\a\ad L_{-m})}\sum\limits_{p=0}^{\infty}
(-\a)^{p}r!\left[\begin{array}{c} i+m \\ p
\end{array}
\right]_m \left[\begin{array}{c} i+(r-p-2)m \\ r
\end{array}
\right]_{m} L_{i+(r-q)m}\nonumber
\\
\!\!\!&=\!\!\!&\sum\limits_{p=0}^{\infty}\sum\limits_{q=0}^{\infty}(-1)^p\a^{q+p}r!\left[\begin{array}{c}
i+m \\ p
\end{array}
\right]_m \left[\begin{array}{c} i+(r-p-2)m \\ r
\end{array}
\right]_{m}\times\nonumber\\
 \!\!\!& \!\!\!& \times\left[\begin{array}{c} i+(r-p+1)m  \\ q
\end{array}
\right]_{m} L_{i+(r-p-q)m}\nonumber\\
\label{to-prove-com}
\!\!\!&=\!\!\!&\sum\limits_{q=0}^{\infty}\sum\limits_{p=0}^{q}(-1)^p\a^{q}r!
\left[\begin{array}{c} i+m \\ p
\end{array}
\right]_m \left[\begin{array}{c} i+(r-p-2)m \\ r
\end{array}
\right]_{m} \times\nonumber\\
 \!\!\!& \!\!\!& \times\left[\begin{array}{c} i+(r-p+1)m \\ q-p
\end{array}
\right]_{m}L_{i+(r-q)m},
\end{eqnarray}
where the last equality follows by first setting $q=q'-p$ and
exchanging summands over $q'$ and $p$ and then replacing $q'$ by
$q$. From the fact that $(\ad Y)^r(L_i)\in\oplus_{q=i-rm}^{i+rm}\F
L_q$, we prove (\ref{3.4'}). Note that
\begin{eqnarray*}
&\!\!\!\!\!\!\!\! &\exp {(\a\ad
L_{-m})}(G_k)=\sum\limits_{p=0}^{\infty}(-\a)^p\left[\begin{array}{c} k+\frac{1}{2}m\\
p
\end{array}
\right]_{m} G_{k-pm},\\
&\!\!\!\!\!\!\!\! &(\ad Y')^r (G_k)=r!\left[\begin{array}{c} k+(r-\frac{3}{2})m\\
r
\end{array}
\right]_{m} G_{k+rm}.
\end{eqnarray*}
We have
\begin{eqnarray}\label{to-prove-com'}
(\ad Y)^r(G_k)
\!\!\!&=\!\!\!&\sum\limits_{q=0}^{\infty}\sum\limits_{p=0}^{q}(-1)^p\a^{q}r!
\left[\begin{array}{c} k+\frac{m}{2} \\ p
\end{array}
\right]_m \left[\begin{array}{c} k+(r-p-\frac32)m \\ r
\end{array}
\right]_{m} \times\nonumber\\
 \!\!\!& \!\!\!& \times\left[\begin{array}{c} k+(r-p+\frac12)m \\ q-p
\end{array}
\right]_{m}G_{k+(r-q)m},
\end{eqnarray}
Then (\ref{3.4''}) follows. \QED \vskip7pt

As a by product of (\ref{to-prove-com}) and (\ref{to-prove-com'}),
we immediately obtain the combinatorial identities (\ref{combin1})
and (\ref{combin2}).

 \begin{lemm}
For any $a \in \F, r, s\in\Z_+$, $i\in\Z$, and $k\in\frac12\Z$, the
following equations hold in $\UU (\LL)$.
\begin{eqnarray}
\label{3.1}&\!\!\!\!\!\!\!\! & L_i X_a ^{<r>} = \sum\limits_{p=0}^r\frac{\a^pr!}{(r-p)!}\left[\begin{array}{c} (p-2)m-i \\
p
\end{array} \right]_m X_{a+p-\frac{i}{m}}^{<r-p>}L_{i-pm},
\\
\label{3.1'}&\!\!\!\!\!\!\!\! & G_k X_a ^{<r>} = \sum\limits_{p=0}^r\frac{\a^pr!}{(r-p)!}\left[\begin{array}{c} (p-\frac32)m-k \\
p
\end{array} \right]_m X_{a+p-\frac{k}{m}}^{<r-p>}G_{k-pm},
\\\label{3.2}&\!\!\!\!\!\!\!\! & L_i X_a ^{[r]}= \sum\limits_{p=0}^r
\frac{\a^pr!}{(r-p)!}\left[\begin{array}{c} (p-2)m-i \\
p
\end{array} \right]_m X_{a-\frac{i}{m}}^{[r-p]}L_{i-pm},
\\\label{3.2'}&\!\!\!\!\!\!\!\! & G_k X_a ^{[r]}= \sum\limits_{p=0}^r
\frac{\a^p r!}{(r-p)!}\left[\begin{array}{c} (p-\frac{3}{2})m-k \\
p
\end{array} \right]_m X_{a-\frac{k}{m}}^{[r-p]}G_{k-pm},
\\
\label{3.3}&\!\!\!\!\!\!\!\! & Y^s X_a ^{<r>} = X_{a-s } ^{<r>}
Y^{s} , \quad\quad\quad Y^s X_a ^{[r]} = X_{a-s } ^{[r]} Y^{s} ,
 \\\label{3.5'}&\!\!\!\!\!\!\!\! &
 L_i Y ^{r} = \sum\limits_{p=0}^r\sum\limits_{q=0}^{2p}
 \frac{(-1)^p\a^qr!}{(r-p)!}a_q(p,i)Y^{r-p}L_{i+(p-q)m},\\\label{3.5''}&\!\!\!\!\!\!\!\! &
 G_k Y ^{r} = \sum\limits_{p=0}^r\sum\limits_{q=0}^{2p}
 \frac{(-1)^p\a^qr!}{(r-p)!}b_q(p,k)Y^{r-p}G_{k+(p-q)m}
.
 \end{eqnarray}
\end{lemm}

\ni{\it Proof.~} We prove (\ref{3.1}) and (\ref{3.1'}) by induction on $r$. The case
of $r=1$ follows from the formula $L_i X =(X-\frac{i}{m}) L_i-\a
(m+i)L_{i-m}$ and $G_k X =(X-\frac{k}{m}) G_k-\a
(\frac{m}{2}+k)G_{k-m}.$ Let
$$
a_p^r=\frac{\a^pr!}{(r-p)!}\left[\begin{array}{c} (p-2)m-i \\
p
\end{array} \right]_m X_{a+p-\frac{i}{m}}^{<r-p>},\quad\quad
b_p^r=\frac{\a^pr!}{(r-p)!}\left[\begin{array}{c} (p-\frac{3}{2})m-k \\
p
\end{array} \right]_m X_{a+p-\frac{k}{m}}^{<r-p>}.
$$
Suppose that (\ref{3.1}) and (\ref{3.1'}) holds for $r$. As for the
case $r+1$, we have
$$\begin{array}{lll}L_i X_a ^{<r+1>} \!\!\!&= \!\!\!&L_i X_a
^{<r>}(X+a+r)\\[4pt]
\!\!\!&= \!\!\!&\sum\limits_{p=0}^r
a_p^r\big((X+a+r+p-\frac{i}{m})L_{i-pm}+\a (pm-m-i)L_{i-(p+1)m}\big)\\[4pt]
\!\!\!&= \!\!\!&\sum\limits_{p=0}^{r+1}\big(a_p^r(X+a+r+p-\frac{i}{m})+\a a_{p-1}^r(pm-2m-i)\big)L_{i-pm}\\[4pt]
\!\!\!&= \!\!\!&\sum\limits_{p=0}^{r+1}
a_p^{r+1}L_{i-pm}\end{array}$$ and
$$\begin{array}{lll}G_k X_a
^{<r+1>} \!\!\!&= \!\!\!&G_k X_a
^{<r>}(X+a+r)\\[4pt]
\!\!\!&= \!\!\!&\sum\limits_{p=0}^r
b_p^r\big((X+a+r+p-\frac{k}{m})G_{k-pm}+\a (pm-\frac12 m-k)G_{k-(p+1)m}\big)\\[4pt]
\!\!\!&= \!\!\!&\sum\limits_{p=0}^{r+1}\big(b_p^r(X+a+r+p-\frac{k}{m})+\a b_{p-1}^r(pm-\frac32 m-k)\big)G_{k-pm}\\[4pt]
\!\!\!&= \!\!\!&\sum\limits_{p=0}^{r+1}
b_p^{r+1}G_{k-pm}.\end{array}$$ Then (\ref{3.1}) and (\ref{3.1'})
hold. Let
$$c_p^r=\frac{\a^pr!}{(r-p)!}\left[\begin{array}{c} (p-2)m-i \\
p
\end{array} \right]_m X_{a-\frac{i}{m}}^{[r-p]},\quad\quad
d_p^r=\frac{\a^pr!}{(r-p)!}\left[\begin{array}{c} (p-\frac32)m-k \\
p
\end{array} \right]_m X_{a-\frac{k}{m}}^{[r-p]}.
$$
From the fact that $$c_p^{r+1}=c_p^r(X+a-r+p-\frac{i}{m})+\a
c_{p-1}^r(pm-2m-i)$$ and $$d_p^{r+1}=d_p^r(X+a-r+p-\frac{k}{m})+\a
d_{p-1}^r(pm-\frac32 m-k),$$ we can prove (\ref{3.2}) and
(\ref{3.2'}) by induction on $r$. The proof of (\ref{3.3}) is
similar.
Using (\ref{3.4'}), together with the fact that
$$L_iY^r=\sum\limits_{p=0}^r(-1)^p\dbinom{r}{p}Y^{r-p}({\rm ad
}Y)^p(L_i),$$ we immediately get (\ref{3.5'}). The proof of
(\ref{3.5''}) is similar with (\ref{3.5'}).  \QED

The following lemma belongs to \cite{G}.

\begin{lemm}\label{lem2} For any $x\in\UU(\LL)$, $a, d \in \F $ and
 $r, s, m\in \Z_+$, we have
 \begin{eqnarray}
\label{2.5}&\!\!\!\!\!\!\!\! &
 x_a ^{<r+ s>} = x_a ^{<r>} x_{a+r} ^{<s>},
\\[6pt]
\label{2.6}&\!\!\!\!\!\!\!\! & x_a ^{[r + s]} = x_a ^{[r]} x_{a-r}
^{[s]},
\\[6pt]
\label{2.7}&\!\!\!\!\!\!\!\! &
 x_a ^{[r]} = x_{a-r+1} ^{<r>},
\\[6pt]
\label{2.8}&\!\!\!\!\!\!\!\! &
 \sum\limits_{r+s = m} \frac{(-1)^s }{r!s!} x_a ^{[r]} x_d ^{<s>} = \dbinom{a-d}{m},
\\[6pt]
\label{2.9}&\!\!\!\!\!\!\!\! &
 \sum\limits_{r+s=m} \frac{(-1)^s }{r!s!} x_a ^{[r]} x_{d-r} ^{[s]}= \dbinom{a-d+m-1}{m} .
 \end{eqnarray}
 \end{lemm}

For $a \in \F$, we set \begin{eqnarray*}&\!\!\!\!\!\!\!\! & \FF_a
=\sum\limits_{r= 0 } ^{\infty } \frac{(-1)^{r}}{r!} X_a ^{[r] }
\otimes Y ^r t^r  , \ \ \ \ F_a = \sum\limits_{r = 0}
^{\infty}\frac{1 }{r! } X_a ^{<r> } \otimes Y^r t^r, \\[4pt]
&\!\!\!\!\!\!\!\! & u_a = \mu \cdot (S_0 \otimes Id ) (F_a ), \ \ \
\ \ \ \ \ \ \  \  v_a = \mu \cdot (Id \otimes S_0 ) (\FF_a
).\end{eqnarray*} Since $S_0 (X_a ^{<r> } ) = (-1 )^r X_{-a } ^{[r]
}$ and $ S_0 (Y ^{r} ) = (-1)^{r} Y^r$, we have
\begin{eqnarray*} u_a = \sum\limits_{r=0} ^{\infty} \frac{(-1)^r
}{r!} X_{-a} ^{[r]} Y ^r t^r ,\ \ \ \  v_a  = \sum\limits_{r=0}
^{\infty} \frac{1}{r!} X_{a} ^{[r]} Y ^r t^r.
\end{eqnarray*}
Denote $\FF = \FF_0 , F = F_0 , u = u_0 , v = v_0 .$ Following the
results in \cite{GZ} and \cite{Y}, $\FF$ is a Drinfel'd twist of the
Witt algebra. The observation of $\WW\subset \LL$ implies that $\FF$
is also a Drinfel'd twist of $\LL$.

\begin{lemm}\adddot
\label{lemma3.2} For $a, d \in \F, $ we have $$\FF_a F_d = 1 \otimes
(1 - Y t )^{a- d}, \ \ \ v_a u_d = (1- Y t)^{-(a + d)}.$$ \noindent
Therefore  $F_a , \FF_a ,  u_a , v_a$ are all invertible, and $\FF_a
^{-1} = F_{a} , u_a ^{-1} = v_{-a}$.
\end{lemm}
\ni{\it Proof.} Using (\ref{2.8}), we have
\begin{eqnarray*}\FF_a F_d\!\!\!&= \!\!\!&(\sum\limits_{r=0} ^{\infty} \frac{(-1)
^r}{r!} X_a ^{[r]} \otimes Y^r t^r ) \cdot (\sum\limits_{s=0}
^{\infty} \frac{1}{s!} X_d ^{<s>} \otimes Y^s t^s )\\\!\!\!&=
\!\!\!&\sum\limits_{r, s =0} ^{\infty} \frac{(-1)^{r}}{r!s!} X_{a}
^{[r]
} X_d ^{<s>} \otimes Y^{r+s} t^{r+s} \\[4pt]
\!\!\!&= \!\!\!&\sum\limits_{p=0} ^{\infty} (-1)^p
\dbinom{a-d}{p} \otimes Y ^{p} t^p \\[4pt]
\!\!\!&= \!\!\!&1 \otimes (1- Y t)^{a-d}.\end{eqnarray*} Using
(\ref{2.9}) and the second formula of (\ref{3.3}), we have
\begin{eqnarray*} v_a u_d \!\!\!&=
\!\!\!&(\sum\limits_{r=0} ^{\infty} \frac{1}{r!} X_{a} ^{[r]} Y ^r
t^r )\cdot(\sum\limits_{s=0} ^{\infty} \frac{(-1)^s
}{s!} X_{-d} ^{[s]} Y ^s t^s) \\[4pt]
\!\!\!&= \!\!\!&\sum\limits_{r, s = 0} ^{\infty}
\frac{1}{r!}\frac{(-1)^{s}}{s!} X_a ^{[r]} Y ^r
X_{-d} ^{[s]} Y ^s t^{r+s} \\[4pt]
\!\!\!&= \!\!\!&\sum\limits_{p=0} ^{\infty}\dbinom{a+d+p-1}{p}
Y ^{p} t^p \ \ \  \\[4pt]
\!\!\!&= \!\!\!& (1-Y t)^{-(a+d)}. \end{eqnarray*} Then this Lemma
follows. \hfill $\Box$

\begin{lemm}\adddot
\label{lemma3.5} For $a \in \F, i \in \Z,$ we have
\begin{eqnarray}
\label{3.7}\!\!\!\!\!\!\!\!\!\!\!\!\!\!\!\!&\!\!\!\!\!\!\!\! &
(L_{i} \otimes 1 )F_a = \sum\limits_{s=0}^{\infty}\a^s\left[\begin{array}{c} (s-2)m-i \\
s
\end{array} \right]_m F_{a-\frac{i}{m}+s}\big(L_{i-sm}\otimes
Y^st^s\big),
\\
\label{3.7'}\!\!\!\!\!\!\!\!\!\!\!\!\!\!\!\!&\!\!\!\!\!\!\!\! &
(G_{k} \otimes 1 )F_a = \sum\limits_{s=0}^{\infty}\a^s\left[\begin{array}{c} (s-\frac32)m-k \\
s
\end{array} \right]_m F_{a-\frac{k}{m}+s}\big(G_{k-sm}\otimes
Y^st^s\big),
\\
\label{3.8}\!\!\!\!\!\!\!\!\!\!\!\!\!\!\!\!&\!\!\!\!\!\!\!\! &(1
\otimes L_{i} )
F_a=\sum\limits_{s=0}^{\infty}(-1)^sF_{a+s}\bigg(\sum\limits_{p=0}^{2s}\a^p
a_p(s,i)X_a^{<s>}\otimes L_{i+(s-p)m}t^{s}\bigg),\\
\label{3.8'}\!\!\!\!\!\!\!\!\!\!\!\!\!\!\!\!&\!\!\!\!\!\!\!\! &(1
\otimes G_{k} )
F_a=\sum\limits_{s=0}^{\infty}(-1)^sF_{a+s}\bigg(\sum\limits_{p=0}^{2s}\a^p
b_p(s,i)X_a^{<s>}\otimes G_{k+(s-p)m}t^{s}\bigg),\\
\!\!\!\!\!\!\!\!\!\!\!\!\!\!\!\!&\!\!\!\!\!\!\!\! &\label{3.9} L_{i}
u_a = u_{a+\frac{i}{m}}\sum\limits_{s=0} ^{\infty}
\sum\limits_{p=0} ^{\infty}\sum\limits_{q=0}^{2p}(-1)^{s}\a^{s+q}\left[\begin{array}{c} (s-2)m-i \\
s
\end{array} \right]_m\times\\\nonumber
\!\!\!\!\!\!\!\!\!\!\!\!\!\!\!\!&\!\!\!\!\!\!\!\! &\phantom{L_{i}
u_a = }\times a_q(p,i-sm)X_{-a-\frac{i}{m}}^{[p]}L_{i+(p-s-q)m}Y
^{s} t^{s+p},\\
\!\!\!\!\!\!\!\!\!\!\!\!\!\!\!\!&\!\!\!\!\!\!\!\! &\label{3.9'}
G_{k} u_a = u_{a+\frac{k}{m}}\sum\limits_{s=0} ^{\infty}
\sum\limits_{p=0} ^{\infty}\sum\limits_{q=0}^{2p}(-1)^{s}\a^{s+q}\left[\begin{array}{c} (s-\frac32)m-k \\
s
\end{array} \right]_m\times\\\nonumber
\!\!\!\!\!\!\!\!\!\!\!\!\!\!\!\!&\!\!\!\!\!\!\!\! &\phantom{G_{k}
u_a = }\times b_q(p,k-sm)X_{-a-\frac{k}{m}}^{[p]}G_{k+(p-s-q)m}Y
^{s} t^{s+p}.
\end{eqnarray}
\end{lemm}
\ni{\it Proof.}  From (\ref{3.1}) and the definition of $F_a$, we
have
\begin{eqnarray*}(L_{i} \otimes 1) F_a
 \!\!\!&= \!\!\!& \sum\limits_{r=0}^{\infty}\frac{1}{r!}L_iX_a^{<r>}\otimes
Y^rt^r\\[4pt]
 \!\!\!&= \!\!\!& \sum\limits_{r=0}^{\infty}\frac{1}{r!}
 \sum\limits_{s=0}^{r}\frac{\a^sr!}{(r-s)!}\left[\begin{array}{c} (s-2)m-i \\
s
\end{array} \right]_m X_{a-\frac{i}{m}+s}^{<r-s>}L_{i-sm}\otimes Y^rt^r
\end{eqnarray*}
\begin{eqnarray*} \!\!\!&= \!\!\!&
\sum\limits_{r=0}^{\infty}\sum\limits_{s=0}^{\infty}\frac{\a^s}{r!}
\left[\begin{array}{c} (s-2)m-i \\
s
\end{array} \right]_m X_{a-\frac{i}{m}+s}^{<r>}L_{i-sm}\otimes Y^{r+s}t^{r+s} \\[4pt]
 \!\!\!&= \!\!\!&\sum\limits_{s=0}^{\infty}\a^s\left[\begin{array}{c} (s-2)m-i \\
s
\end{array} \right]_m\bigg(\sum\limits_{r=0}^{\infty}\frac{1}{r!}X_{a-\frac{i}{m}+s}^{<r>}\otimes
Y^{r}t^{r}\bigg)
\big(L_{i-sm}\otimes Y^st^s\big)\\[4pt]
\!\!\!&=\!\!\!&\sum\limits_{s=0}^{\infty}\a^s\left[\begin{array}{c} (s-2)m-i \\
s
\end{array} \right]_m F_{a-\frac{i}{m}+s}\big(L_{i-sm}\otimes
Y^st^s\big).
\end{eqnarray*}
This proves (\ref{3.7}). Similarly, (\ref{3.7'}) follows from
(\ref{3.1'}). For (\ref{3.8}), using (\ref{2.5}) and (\ref{3.5'}),
we have
\begin{eqnarray*} (1 \otimes L_i ) F_a \!\!\!&= \!\!\!&
\sum\limits_{r=0}^{\infty}\frac{1}{r!}X_a^{<r>}\otimes
L_iY^rt^r \\[4pt]
\!\!\!&=
\!\!\!&\sum\limits_{r=0}^{\infty}\sum\limits_{s=0}^{\infty}\frac{(-1)^s}{r!}\sum\limits_{p=0}^{2s}
\a^p a_p(s,i)X_a^{<r+s>}\otimes Y^{r}L_{i+(s-p)m}t^{r+s}\\[4pt]
\!\!\!&=\!\!\!&\sum\limits_{s=0}^{\infty}(-1)^sF_{a+s}\bigg(\sum\limits_{p=0}^{2s}\a^p
a_p(s,i)X_a^{<s>}\otimes L_{i+(s-p)m}t^{s}\bigg),
\end{eqnarray*}
So (\ref{3.8}) is right. (\ref{3.8'}) is obtained from (\ref{2.5})
and (\ref{3.5''}). Now we prove (\ref{3.9}). From (\ref{2.6}),
(\ref{3.2}), (\ref{3.3}) and (\ref{3.5'}), we
 get
\begin{eqnarray*} L_{i} u_a \!\!\!& =\!\!\!&  \sum\limits_{r=0} ^{\infty}
\frac{(-1)^r}{r!} L_{i}X_{-a} ^{[r]}Y ^rt^r \\
\!\!\!&= \!\!\!& \sum\limits_{r=0} ^{\infty} \frac{(-1)^r}{r!}
\sum\limits_{s=0} ^{r}
\frac{\a^sr!}{(r-s)!}\left[\begin{array}{c} (s-2)m-i \\
s
\end{array} \right]_mX_{-a-\frac{i}{m}} ^{[r-s]}L_{i-sm} Y ^r t^r \\
\!\!\!&= \!\!\!& \sum\limits_{r=0} ^{\infty} \sum\limits_{s=0}
^{\infty}
\frac{(-1)^{r+s}\a^s}{r!}\left[\begin{array}{c} (s-2)m-i \\
s
\end{array} \right]_mX_{-a-\frac{i}{m}} ^{[r]}L_{i-sm} Y ^{r+s} t^{r+s} \\
\!\!\!&= \!\!\!& \sum\limits_{r=0} ^{\infty} \sum\limits_{s=0}
^{\infty}
\frac{(-1)^{r+s}\a^s}{r!}\left[\begin{array}{c} (s-2)m-i \\
s
\end{array} \right]_mX_{-a-\frac{i}{m}} ^{[r]}\times\\
\!\!\!& \!\!\!&\times \sum\limits_{p=0}^r\sum\limits_{q=0}^{2p}
\frac{(-1)^p\a^qr!}{(r-p)!}a_q(p,i-sm)Y^{r-p}L_{i+(p-s-q)m}Y ^{s} t^{r+s} \\
\!\!\!&= \!\!\!& \sum\limits_{s=0} ^{\infty} \sum\limits_{p=0}
^{\infty}\bigg(\sum\limits_{r=0}^{\infty}
\frac{(-1)^{r}}{r!}X_{-a-\frac{i}{m}} ^{[r]}Y^rt^r\bigg)(-1)^{s}\left[\begin{array}{c} (s-2)m-i \\
s
\end{array} \right]_mX_{-a-\frac{i}{m}}^{[p]}\times\\
\!\!\!& \!\!\!&\times
\sum\limits_{q=0}^{2p}\a^{s+q}a_q(p,i-sm)L_{i+(p-s-q)m}Y ^{s}
t^{s+p}
\end{eqnarray*}
\begin{eqnarray*}\!\!\!&= \!\!\!& u_{a+\frac{i}{m}}\sum\limits_{s=0} ^{\infty}
\sum\limits_{p=0} ^{\infty}\sum\limits_{q=0}^{2p}(-1)^{s}\a^{s+q}\left[\begin{array}{c} (s-2)m-i \\
s
\end{array} \right]_m
a_q(p,i-sm)X_{-a-\frac{i}{m}}^{[p]}L_{i+(p-s-q)m}Y ^{s} t^{s+p}.
\end{eqnarray*}
Thus (\ref{3.9}) holds. Similarly,  (\ref{3.9'}) can be proved from
 (\ref{3.2'}), (\ref{3.3}),(\ref{3.5''}), and (\ref{2.6}).\QED

Now we give the proof of Theorem \ref{main}.
 \vskip10pt
\ni{\it Proof of Theorem \ref{main}.} Sice $\FF$ is a Drinfel'd
twist, according to Lemma \ref{lemma2.3}, we only need to determine
the action of $\D$ and $S$ on $L_i, G_k \in \LL$, $ i\in\Z,
k\in\frac12\Z$. From (\ref{3.7}), (\ref{3.8}), and Lemma
\ref{lemma3.2}, we have
\begin{eqnarray*}\D(L_{i} )\!\!\!&= \!\!\!& \FF \D_0 (L_i) \FF^{-1} = \FF (L_i \otimes 1) \FF^{-1} + \FF (1 \otimes L_i)
\FF^{-1} =\FF (L_i \otimes 1) F + \FF (1 \otimes L_i)
F \\
\!\!\!&= \!\!\!&\FF\sum\limits_{r=0}^{\infty}\a^r\left[\begin{array}{c} (r-2)m-i \\
r
\end{array} \right]_m F_{-\frac{i}{m}+r} \big(L_{i-rm}  \otimes Y^rt^r) \\
\!\!\!& \!\!\!& +\FF
\sum\limits_{r=0}^{\infty}(-1)^rF_{r}\bigg(\sum\limits_{s=0}^{2r}
\a^sa_s(r,i)X^{<r>}\otimes L_{i+(r-s)m}t^{r}\bigg)\\
\!\!\!&= \!\!\!& \sum\limits_{r=0}^{\infty}\a^r\left[\begin{array}{c} (r-2)m-i \\
r
\end{array} \right]_m\big(1\otimes (1-Yt)^{\frac{i}{m}-r}\big)(L_{i-rm}\otimes Y^rt^r)\\
\!\!\!& \!\!\!& +\sum\limits_{r=0}^{\infty}(-1)^r\big(1\otimes
(1-Yt)^{-r}\big)\bigg(\sum\limits_{s=0}^{2r}\a^s
a_s(r,i)X^{<r>}\otimes L_{i+(r-s)m}t^r\bigg)\\
\!\!\!& =\!\!\!&\sum\limits_{r=0}^{\infty}\a^r\left[\begin{array}{c} (r-2)m-i \\
r
\end{array} \right]_m\big(L_{i-rm}\otimes (1-Yt)^{\frac{i}{m}-r}Y^rt^r\big)\\
\!\!\!& \!\!\!&
+\sum\limits_{r=0}^{\infty}(-1)^r\sum\limits_{s=0}^{2r}\a^s
a_s(r,i)\big(X^{<r>}\otimes (1-Yt)^{-r}L_{i+(r-s)m}t^r\big).
\end{eqnarray*}
From (\ref{3.7'}), (\ref{3.8'}), and Lemma \ref{lemma3.2}, we have
\begin{eqnarray*}\D(G_{k} )\!\!\!&= \!\!\!& \FF \D_0 (G_k) \FF^{-1}
= \FF (G_k \otimes 1) \FF^{-1} + \FF (1 \otimes G_k) \FF^{-1} =\FF
(G_k \otimes 1) F + \FF (1 \otimes G_k)
F \\
\!\!\!&= \!\!\!&\FF\sum\limits_{r=0}^{\infty}\a^r\left[\begin{array}{c} (r-\frac32)m-k \\
r
\end{array} \right]_m F_{-\frac{k}{m}+r} \big(G_{k-rm}  \otimes Y^rt^r) \\
\!\!\!& \!\!\!& +\FF
\sum\limits_{r=0}^{\infty}(-1)^rF_{r}\bigg(\sum\limits_{s=0}^{2r}
\a^sb_s(r,k)X^{<r>}\otimes G_{k+(r-s)m}t^{r}\bigg)\\
\!\!\!&= \!\!\!& \sum\limits_{r=0}^{\infty}\a^r\left[\begin{array}{c} (r-\frac32)m-k \\
r
\end{array} \right]_m\big(1\otimes (1-Yt)^{\frac{k}{m}-r}\big)(G_{k-rm}\otimes Y^rt^r)\\
\!\!\!& \!\!\!& +\sum\limits_{r=0}^{\infty}(-1)^r\big(1\otimes
(1-Yt)^{-r}\big)\bigg(\sum\limits_{s=0}^{2r}\a^s
b_s(r,k)X^{<r>}\otimes G_{k+(r-s)m}t^r\bigg)
\end{eqnarray*}
\begin{eqnarray*}\!\!\!& =\!\!\!&\sum\limits_{r=0}^{\infty}\a^r\left[\begin{array}{c} (r-\frac32)m-k \\
r
\end{array} \right]_m\big(G_{k-rm}\otimes (1-Yt)^{\frac{k}{m}-r}Y^rt^r\big)\\
\!\!\!& \!\!\!&
+\sum\limits_{r=0}^{\infty}(-1)^r\sum\limits_{s=0}^{2r}\a^s
b_s(r,k)\big(X^{<r>}\otimes (1-Yt)^{-r}G_{k+(r-s)m}t^r\big).
\end{eqnarray*}

Using (\ref{3.9}) and Lemmas \ref{lemma3.2}, we have
\begin{eqnarray*}S(L_i) \!\!\!&= \!\!\!& u^{-1} S_{0} (L_i) u
=- v L_i u\\
\!\!\!&= \!\!\!& -v u_{\frac{i}{m}}\sum\limits_{r=0} ^{\infty}
\sum\limits_{p=0} ^{\infty}\sum\limits_{q=0}^{2p}(-1)^{r}\a^{r+q}\left[\begin{array}{c} (r-2)m-i \\
r
\end{array} \right]_m
a_q(p,i-rm)X_{-\frac{i}{m}}^{[p]}L_{i+(p-r-q)m}Y ^{r} t^{r+p}\\
\!\!\!&= \!\!\!& -(1-Yt)^{-\frac{i}{m}}\sum\limits_{r=0} ^{\infty}
\sum\limits_{p=0} ^{\infty}\sum\limits_{q=0}^{2p}(-1)^{r}\a^{r+q}\left[\begin{array}{c} (r-2)m-i \\
r
\end{array} \right]_m\times\\
\!\!\!& \!\!\!&\times
a_q(p,i-rm)X_{-\frac{i}{m}}^{[p]}L_{i+(p-r-q)m}Y ^{r} t^{r+p}.
\end{eqnarray*}

Using (\ref{3.9'}) and Lemmas \ref{lemma3.2}, we have
\begin{eqnarray*}S(G_k) \!\!\!&= \!\!\!& u^{-1} S_{0} (G_k) u
=- v G_k u\\
\!\!\!&= \!\!\!& -v u_{\frac{k}{m}}\sum\limits_{r=0} ^{\infty}
\sum\limits_{p=0} ^{\infty}\sum\limits_{q=0}^{2p}(-1)^{r}\a^{r+q}\left[\begin{array}{c} (r-\frac32)m-k \\
r
\end{array} \right]_m
b_q(p,k-rm)X_{-\frac{k}{m}}^{[p]}G_{k+(p-r-q)m}Y ^{r} t^{r+p}\\
\!\!\!&= \!\!\!& -(1-Yt)^{-\frac{k}{m}}\sum\limits_{r=0} ^{\infty}
\sum\limits_{p=0} ^{\infty}\sum\limits_{q=0}^{2p}(-1)^{r}\a^{r+q}\left[\begin{array}{c} (r-\frac32)m-k \\
r
\end{array} \right]_m\times\\
\!\!\!& \!\!\!&\times
b_q(p,k-rm)X_{-\frac{k}{m}}^{[p]}G_{k+(p-r-q)m}Y ^{r} t^{r+p}.
\end{eqnarray*}
The proof is completed.\qed

\end{document}